\def\LS{{\cal LS}}
 \magnification = \magstep1
\baselineskip =13pt
\overfullrule =0pt

\centerline{\bf INJECTIVE RESOLUTIONS OF $BG$  AND DERIVED}
\centerline{\bf  MODULI SPACES OF LOCAL SYSTEMS}
\vskip .7cm

\centerline{\bf M. Kapranov}

\vskip 1cm

Let $X$ be a finite connected CW-complex, $x_0\in X$ be a point,
and $G$ be an affine algebraic group
over {\bf C}. A $G$-local system on $X$ is just a locally constant sheaf
of $G$-torsors. Let $\LS_G(X, x_0)$ denote the set
of isomorphism classes of $G$-local systems trivialized over $x_0$.
For such a local system $E$ let $[E]$ denote the corresponding point
of $\LS_G(X, x_0)$. 
This set is naturally an algebraic variety; it is just the
variety of all homomorphisms $\pi_1(X, x_0)\to G$. It is acted
upon
by $G$, and the quotient $\LS_G(X) = \LS_G(X, x_0)/G$
is  the set of isomorphism classes of local systems. 
Since the $G$-action may not be free, $\LS_G(X)$ may not exist
as an algebraic variety, but is well defined as an algebraic stack.
The first order deformation theory gives an identification
$T_{[E]}\LS_G(X) = H^1(X, {\rm ad}(E))$ for any local system.
The stack $\LS_G(X)$ and the variety $\LS_G(X, x_0)$ may be not smooth:
the jumping of the dimension of the tangent space is made possible
by the corresponding jumping of the dimensions of the higher cohomology
of ${\rm ad}(E)$.
 
\vskip .2cm

This situation is typical for many other problems of deformation
theory (e.g., the moduli stack of algebraic vector bundles on a variety $M$
is smooth when $M$ is a  curve, but not when $\dim(M)>1$). It has been
proposed by several people  (among them Deligne, Drinfeld and Kontsevich,
see [K]) that one could overcome this difficulty by systematically
working in the derived category, i.e., constructing a kind of non-Abelian
derived functor of the moduli space functor. The appropriate language
for such derived functors (over a field of characteristic 0)
is that of dg-schemes, i.e., schemes together with a sheaf of
(negatively graded) differential graded algebras, cf. [Q], [Mun].
From this
point of view the reason that the moduli space is singular is that
we disregard the higher cohomology by artificially truncating the
``true" derived moduli space, which should be the right object
to consider in all geometric studies.

\vskip .2cm

It is not, however, exactly straightforward to construct
such derived moduli spaces, and the purpose of the present note is to do so
in the simplest case, that of local systems. 
Our construction is based on the observation that the moduli space
space can be represented as $[N, BG]$,
the set of simplicial 
homotopy classes of simplicial maps from a \v Cech nerve $N$
of $X$ to the simplicial classifying space of $G$. 
To extend this into derived category, we construct
an appropriate ``injective" resolution $RBG$ of $BG$. 
It may seem surprising that
the space such as $BG$ needs an additional injective resolution,
since, regarded as a simplicial set, it is fibrant. The reason is that
we are interested in a geometric and not topological concept
of  fibrations for (dg)-schemes.

To keep the paper short, we avoided going into foundational matters related
to dg-schemes; for instance, the derived moduli spaces should really
be ``dg-stacks", but we consider the rigidified situation
when the stack structure is not necessary. We also did not  consider
in any detail the (rudiments of a) closed model structure on
the category of  dg-schemes with the role of fibrations played by smooth 
maps etc. 

\vskip .2cm

\vskip .2cm

I am grateful to M. Kontsevich who convinced me that the dg-point of view
is more flexible than the more exotic approaches to circumventing the
singularity of moduli spaces I was experimenting with. 
This research was supported by an NSF grant and an A.P.Sloan Fellowship.

\vskip 2cm

\centerline {\bf \S 1. Injective resolutions of $BG$.}

\vskip 1cm

\noindent {\bf (1.1) Dg-schemes.} We work everywhere over the field {\bf C}
of complex numbers. By a complex (or dg-vector space) we mean a cochain
complex, i.e., a graded vector space with a differential of degree $+1$. 
By a dg-algebra we always mean a commutative
${\bf Z}_{\leq 0}$-graded differential algebra $A$. Note that in such
an algebra the differential $d$ is $A^0$-linear, so its
cohomology forms a graded $A^0$-module $H^\bullet(A)$.
A quasi-isomorphism of dg-algebras is a morphism inducing an isomorphism
in the cohomology.  

It is also convenient
to use the dual geometric language and to speak about dg-schemes.
By definition, a dg-scheme $M$ is an ordinary scheme $M^0$ equipped
with a quasi-coherent sheaf ${\cal O}_{M}^\bullet$ of dg-algebras
on $M^0$. Thus affine dg-schemes (those with $M^0$ affine)
are in anti-equivalence with dg-algebras. For a dg-scheme $M$
we have the scheme $\pi_0(M) = {\rm Spec}(\underline{H}^0({\cal O}_M))$
which is a closed subscheme in $M^0$. It is possible to
view dg-schemes as superschemes in the sense of [Man] equipped with
additional structure as follows. 

\proclaim (1.1.1) Proposition. Let $Y$ be a scheme. Then the category of 
dg-schemes
$M$ with $M^0=Y$ is equivalent to the category of superschemes $\tilde M\to Y$
affine over $Y$  and equipped with the following additional 
structures:\hfill\break
(1) A section $i: S\to\tilde M$:\hfill\break
(2) An (algebraic) action of the multiplicative semigroup  $({\bf C}, \times)$
 on $\tilde M$ whose fixed point subscheme
is $i(Y)$;\hfill\break\
(3) An odd vector field $d$ on $\tilde M$ satisfying $\{d,d\}=0$ and having 
degree 1,
i.e., such that $[d,L]=L$ where $L$ is the vector field generating the 
action of ${\bf G}_m\i {\bf C}$.

The proof is obvious and left to the reader: the action of ${\bf G}_m$
gives a {\bf Z}-grading, the fact that it is situated in degrees $\leq 0$
is encoded by saying that the action extends to an action of {\bf C} etc. 

Because of this proposition we can easily reduce several foundational
questions regarding dg-schemes (e.g., the properties of the sheaves
of differentials and derivations) to those about superschemes,
which have been treated in [Man].

\vskip .3cm

\noindent {\bf (1.2) Tangent spaces.} 
We will say that a dg-scheme $M$ is smooth (or is a dg-manifold), if $M^0$ is
a smooth manifold, and locally on the Zariski
topology of $M^0$ the sheaf ${\cal O}_M^\bullet$ is free
as a sheaf of graded commutative algebras, with finitely many generators
in each degree. 

Given a dg-manifold $M$ and a
{\bf C}-point $x\in \pi_0(M)$, we have the tangent
 dg-space (complex)
$T_x^\bullet M$. It is defined, as usual, as the graded vector space
of {\bf C}-valued derivations.
 This is a complex of vector
spaces concentrated in degrees $\geq 0$.

Any morphism
$f: M\to N$ of dg-manifolds gives rise to a morphism of complexes
$d_xf: T_x^\bullet M\to T_{f(x)}^\bullet N$.  
It is suggestive
to use the topological notation for the cohomology
of the tangent complex:
$$\pi_{-i}(M,x):= H^{i}(T^\bullet_x M), \quad i\geq 0. \leqno (1.2.1)$$
This notation is justified by the following fact.

\proclaim (1.2.2) Proposition. For any smooth dg-scheme
 $M$ and any its {\bf C}-point
$x$ there are natural bilinear maps (``Whitehead products")
$$\pi_i(M,x)\otimes \pi_j(M, x) \to \pi_{i+j-1}(M,x)$$
which makes $\pi_{\bullet +1}(M,x)$ into a graded Lie algebra. For any morphism
$f: M\to N$ of dg-manifolds the induced morphism $\pi_\bullet(M,x)
\to \pi_\bullet(N, f(x))$ is a homomorphism of graded Lie algebras. 

\noindent {\sl Proof:} Recall the concept of a weak Lie algebra
(or a shLie algebra [St]). This is a  graded vector space {\bf g}
together with a continuous differential $D$ in $\hat S({\bf g}^*[-1])$,
the completed symmetric algebra of the shifted dual vector space. 
By restricting $D$ to the degree 1 part, namely ${\bf g}^*$
and dualizing the graded components of this restriction, one gets
antisymmetric $n$-linear brackets of degree $2-n$
$$\lambda_n: {\bf g}^{\otimes n}\to {\bf g},\,\,
 x_1\otimes ... \otimes x_n \mapsto [x_1, ..., x_n]_n, \,\, n\geq 1,$$
In particular, $d=\lambda_1$ is a differential in {\bf g}, while
$\lambda_2$ satisfies the jacobi identity up to $d$-boundaries, so
that $H^\bullet_d({\bf g})$ is a graded Lie algebra.

Now,  $M$ being a dg-manifold, the completion
$\hat{\cal O}_{M,x}^\bullet$ 
is isomorphic, as a graded algebra, to $\hat S(W^\bullet)$,
the completion of a free graded algebra generated
by some  graded vector space $W^\bullet$. Let $V^\bullet = W^*$ be the dual 
graded
space. 
Such an isomorphism $\phi$ is just a formal coordinate system in $M$ near $x$.
Given such $\phi$,  its differential identifies $V^\bullet$ with $T_x^\bullet 
M$.
So $\phi$ identifies the graded algebra $\hat S(T_x^*M)$ with the dg-algebra
$\hat{\cal O}_{M,x}^\bullet$, and thus we get by pullback a differential $D$ on
$\hat S(T_x^*M)$ which of course satisfies $D^2=0$. 
This means that $T^\bullet_x M[-1]$ becomes equipped with the structure of a  
weak Lie algebra,
so its cohomology is a graded Lie algebra. If we choose a different isomorphism
$\phi': \hat S(V^*)\to \hat{\cal O}_{M,x}^\bullet$ but with the same 
differential
at 0, then we get what is is known as a weakly isomorphic weak Lie algebra,
so that the Lie algebra structure on the cohomology will be the same. 

\vskip .2cm

The following fact can be seen as an analog of the 
Whitehead theorem in topology.

\proclaim (1.2.3) Proposition. Let $f: M\to N$ be a morhism of dg-manifolds.
Then the following conditions are equivalent:\hfill\break
(a) $f$ is a quiasiisomorphism.\hfill\break
(b) The morphism of schemes $\pi_0(f); \pi_0(M)\to \pi_0(N)$ is an isomorphism,
and for any {\bf C}-point $x$ of $M$ the differential $d_xf$ induces
an isomorphism $\pi_i(M,x)\to \pi_i (N, f(x))$ for all $i\leq 0$. 

\noindent {\sl Proof:} It is enough to prove that for any $x\in\pi_0(M)$
the map $f$ induces the completed local dg-algebras
$\hat f^*\hat{\cal O}^\bullet_{N, f(x)}\to \hat{\cal O}^\bullet_{M,x}$,
is a quasiisomorphism.
For that, notice that $\hat{\cal O}^\bullet_{M,x}$ has a filtration
whose quotients are the symmetric powers of the cotangent dg-space
$T^*_xM$. So if $f$ gives a quasiisomorphism of tangent
dg-spaces, we find that $\hat f^*$ induces quasiisomorphisms on
the quotients of the natural filtrations. So the proof is
accomplished by invoking a spectral sequence argument, which
is legitimate (i.e., the spectral sequences converge) because
the dg-algebras in question are ${\bf Z}_{\leq 0}$-graded.

 \vskip .3cm

\noindent {\bf (1.3)  Twisted tensor products and fibrations.} A dg-algebra
$C$ is called a twisted tensor product of dg-algebras $A$ and $B$, if
$C\simeq A\otimes B$ as a graded algebra, and with respect to this
identification,
$$d_C = d_A\otimes 1 + 1\otimes d_B + \sum_{i\geq 2} d_i, \quad
{\rm deg}(d_i) = (1-i, i). \leqno (1.3.1)$$
In this case the natural embedding 
$$A\hookrightarrow C, \quad a\mapsto a\otimes 1, \leqno (1.3.2)$$
is a morphism of dg-algebras. See [May]. Note, in particular, the concept
of a quasi-free dg-algebra $C$ over $A$. This just means that $C$
is a twisted tensor product of $A$ and a free dg-algebra. 

\proclaim (1.3.3) Definition. Let $p: M\to N$ be a morphism of affine 
dg-schemes,
$M={\rm Spec}(C), N={\rm Spec}(A)$. We will say that $p$ is a fibration
with fiber $F={\rm Spec}(B)$, if $C$ is isomorphic to a twisted tensor
product of $A$ and $B$ in such a way that the homomorphism $p^*: A\to C$
becomes the canonical embedding (1.3.2). 

It is clear that for a fibration we have a spectral sequence
$$E_2 = H^\bullet(A)\otimes H^\bullet(B) \Rightarrow H^\bullet(C).
\leqno (1.3.4)$$

\vskip .2cm

\noindent {\bf (1.4) Resolutions.} Let $A\to R$ be any morphism of
  dg-algebras.
Then it is standard, see, e.g., [Mun] how to replace $R$ by a 
quasi-isomorphic 
dg-algebra $\cal R$ which is quasi-free over $A$. 
We construct $\cal R$ as the union
of an increasing sequence ${\cal R}_n, n\geq 0$ of sub-dg-algebras.
 We first take a subspace of algebra
generators in $H^0(R)$, and  lift this space
to  a subspace $V^0\i R^0$. We define ${\cal R}_0$ to be the free
$A$-algebra generated by $V^0$ (where the differential is set to
vanish on $V^0$). Thus we have a morphism of dg-algebras
$d_0: {\cal R}_0\to R$ surjective on $H^0$. 
Then we take a space of generators of the ideal
${\rm Ker}(d_0)$ and denote by $V^{-1}$ the vector space freely
spanned by these generators. Then we have a natural morphism
$d_{-1}: V^{-1}\to {\cal R}_0$ which gives rise to a dg-algebra
structure on the free ${\cal R}_0$-algebra generated by $V^{-1}$. 
 Denote the dg-algebra thus obtained by ${\cal R}_1$.
By construction, $H^0({\cal R}_1) = R$. As the next step,  we take
a space $V^{-2}$ of generators of the $H^0({\cal R}_1)$-module 
$H^{-1}({\cal R}_1)$, lift $V^{-2}$ to a subspace of cocycles
and define in this way
a morphism of graded vector spaces $d_{-2}: V^{-2}\to {\cal R}_1$.
 This gives
a quasi-free dg-algebra ${\cal R}_2$. Continuing in this
way, we  inductively construct ${\cal R}_n$ so as to kill the
$(-n+1)$st cohomology of ${\rm Ker}({\cal R}_{n-1}\to R$ 
and not to affect the
$j$th cohomology, $j> -n+1$. 

Let us summarize the well known properties of this construction.

\proclaim (1.4.1) Proposition. (a) Any two quasi-free resolutions
constructed in this way, are quasiisomorphic.\hfill\break
(b) If $A, R$ have only finitely many generators in each degree,
then we can choose $\cal R$ with the same property.\hfill\break
(c) If, in the situation of (b),  $H$ is a reductive algebraic group
 acting on $A,R$ so that
the morphism $A\to R$ is equivariant, then it is possible to
choose $\cal R$ so as to possess a $G$-action compatible with the
maps and to have finitely many generators in each degree. 

To see part (c), just notice that it is possible to take
the spaces of generators on each step of construction
to be $H$-invariant and finite-dimensional.

\vskip .3cm

\noindent {\bf (1.5) Simplicial objects and classifying spaces.}
We will be using simplicial objects in the categories of sets, schemes
and dg-schemes. See [F] [May] for general background.
 By $\Delta[n]$ we denote the standard $n$-simplex regarded
as a simplicial set. If $I$ is any set (scheme), then the collection
of Cartesian powers $(I^n)_{n\geq 0}$ forms a simplicial set (scheme)
which we denote $\Delta(I)$ and call the unoriented $I$-simplex. 

\vskip .2cm

Let $\cal C$ be any category. Its nerve
(or classifying space) $B_\bullet {\cal C}$ is the simplicial
set whose $n$-simplices are ``commutative $n$-simplices'' in $\cal C$,
i.e., diagrams consisting of $n+1$ objects $A_0, ..., A_n$ and morphisms
$g_{ij}: A_i\to A_j$ satisfying the
conditions
$$g_{jk}g_{ij}=g_{ik},  \quad  i<j<k. \leqno (1.5.1)$$
Let also  $\tilde B_n{\cal C}$ be the set
consisting of {\it all},  not necessarily commutative,
simplex-shaped diagrams in $\cal C$. In other words, an element of
$\tilde B_n{\cal C}$ is an arbitrary collection of objects $A_0, ..., A_n$
and morphisms $g_{ij}: A_i\to A_j$. It is clear that these sets unite
into a simplicial set $\tilde B_\bullet C$ containing $B_\bullet C$. 

A group $G$ over {\bf C} can be considered as a category
with one object and the set of automorphisms $G$, so $B_\bullet G$ is defined.
If $G$ is an affine algebraic group over ${\bf C}$, then $B_\bullet G$
and $\tilde B_\bullet G$ are simplicial schemes.  Moreover, $B_\bullet G
\to \tilde B_\bullet G$ is a closed embedding, given by the equations
(1.5.1). In other words, for every $n$ we have a surjection
of algebras of functions
$${\bf C}[\tilde B_nG] = \bigotimes_{0\leq i<j\leq n} 
{\bf C}[G] \longrightarrow {\bf C}[B_nG].\leqno (1.5.2)$$
Notice also that the group $G^{n+1}$ acts on $\tilde B_nG$
by ``gauge transformations":
$$(g_0, ..., g_n): (g_{ij})\mapsto (g_jg_{ij}g_i^{-1}), \leqno (1.5.2)$$
and this action preserved $B_nG$. The simplicial scheme $\Delta(G)$
formed by the $G^n$ is actually a simplicial algebraic group, acting
on the simplicial scheme $\tilde BG$ and preserving $BG$.

\vskip .3cm

\noindent {\bf (1.6) Injective simplicial dg-schemes.} 
A simplicial dg-scheme $X_\bullet$ is called affine if each $X_i$
is affine. Such a scheme is the same as a cosimplicial dg-algebra. 

If $S$ is a simplicial set and $X$ is a simplicial dg-scheme,
then we have a dg-scheme $\underline{\rm Hom}(S,X)$, which
is affine when $X$ is affine. 

\proclaim (1.6.1) Definition. An affine dg-scheme $X_\bullet$ is called
injective, if for any cofibration (i.e., embedding) $S'\i S$
of simplicial sets the induced map of dg-schemes
$\underline{\rm Hom}(S,X)\to \underline {\rm Hom}(S',X)$
is smooth and is a fibration with
fiber $\underline{\rm Hom}(S/S', X)$. 

So this definition is a direct analog of the concept of an injective
object in an Abelian category.

\vskip .1cm

\noindent {\bf (1.6.2) Example.} The  simplicial scheme $\tilde B_\bullet G$
(with trivial dg-structure) is injective, but its simplicial subscheme
$B_\bullet G$ is not.

More generally, we have the following fact, whose proof is achieved by
using induction over the simplices. 

\proclaim (1.6.3) Proposition. A simplicial dg-scheme $X$ is injective
if anf only if $X_0$ is smooth and for any $n\geq 1$ the morphism
$$X_n = \underline{\rm Hom}(\Delta [n], X)\to \underline{\rm Hom}
(\partial\Delta[n], X)$$
is a fibration with smooth fiber.

\vskip .2cm

\proclaim (1.7) Theorem. Let $G$ be any reductive group.
One can replace $BG$ by a quasiisomorphic 
(dimension by dimension) injective simplicial dg-manifold $RBG$ with
an action of the simplicial group $\Delta(G)$. This replacement is
canonical up to an equivariant
 quaisiisomorphism of simplicial dg-manifolds. 

\noindent {\sl Proof:} 
We set $RB_mG = B_mG$ for $m=0,1$, then take for ${\bf C}[RB_2G]$
any free dg-resolution of ${\bf C}[B_2G]$ as an algebra over
${\bf C}[\tilde B_2G]$, and then continue inductively as follows.
Suppose we already constructed dg-schemes $RB_nG, n\leq m$
and face morphisms $\partial_i$ satisfying the simplicial identities
(i.e., suppose we constructed the $m$-th skeleton $RB_{\leq m}G$).
Then the dg-scheme $\underline{\rm Hom}(\partial\Delta[m+1], RB_{\leq m}G)$
is defined (because $\partial\Delta[m+1]$ is $m$-dimensional). 
Let $A(m+1)$ be its dg-algebra of functions. We have a natural dg-algebra
morphism $A(m+1)\to {\bf C}[B_{m+1}G]$ (the latter algebra
has, of course, trivial dg-structure). Further, there is a natural action
of $G^{m+2}$ on $A(m+1)$ so that the morphism becomes equivariant.
 Define ${\bf C}[RB_{m+1}G]$
to be an equivariant quasi-free dg-resolution of ${\bf C}[B_{m+1}G]$ as an
$A(m+1)$-algebra. Continuing in this way, we get a required dg-resolution
of the entire $BG$.

\vskip .2cm

\noindent {\bf (1.8) Explicit resolution for $G=GL(r)$.}
In  the case $G=GL(r)$ it is possible to write down a canonical
injective resolution by analyzing the syzygies among the equations
(1.5.1).

 We set $G=GL(r)$. The scheme $\tilde B_nG$ is just the product $\prod_{0\leq 
i<j\leq n} G$
of $n(n+1)/2$ copies of $G$. We will denote the $(i,j)$th copy by $G_{ij}$. 
The ring ${\bf C}[G]$  of functions on $G$ is generated by the matrix elements
of one matrix-valued function $g$ which is required to satisfy $\det(g)\neq 0$.
The ring of functions on $\tilde B_n G$  i.e., $\bigotimes_{0\leq i<j\leq n} 
{\bf C}[G_{ij}]$
 is therefore generated by the matrix
elements of $n(n+1)/2$ matrix-valued functions which we denote by $g_{ij}$, 
$0\leq i<j\leq n$.
The scheme $B_nG$ is described inside $\tilde B_nG$ by $n+1\choose 3$ 
equations
(1.5.1). We now define a dg-algebra ${\bf C}[RB_nG]$ over
 ${\bf C}[\tilde B_nG]$ to be generated
by the matrix elements of the $(r\times r)$-matrix functions
 $g_{i_0...i_p}$, $0\leq i_0 < ... < i_p\leq n$, with the degree of
(each matrix element of) $g_{i_0, ..., i_p}$ equal to $1-p$, and the 
differential $d$
(of degree $+1$) defined on the generators by
$$d(g_{i_0...i_p}) = 
\sum_{\nu=1}^{p-1} (-1)^\nu \bigl(g_{i_0, ..., \hat i_\nu, ..., i_p} - 
g_{i_\nu ...i_p} g_{i_0 ... i_\nu}
\bigr). \leqno (1.8.1)$$
Note in particular that for $p=2$ we get
$$d(g_{ijk}) = g_{jk}g_{ij}-g_{ik},\leqno (1.8.2)$$
so the image of the last differential is the ideal in ${\bf C}[\tilde B_nG]$
generated by the equations (1.5.1).

One verifies immediately that the condition $d^2=0$ is satisfied on the 
generators
(and hence on the entire algebra). So ${\bf C}[RB_nG]$ is a free dg-algebra 
over
${\bf C}[\tilde B_nG]$. We denote by $RB_nG$ the affine dg-scheme whose ring
of functions is ${\bf C}[RB_nG]$. 

\proclaim (1.8.3) Theorem. (a) The dg-algebra ${\bf C}[RB_nG]$ is 
quasiisomorphic to
${\bf C}[B_nG]$, i.e., it provides a free resolution of the equations (1.5.1). 
\hfill\break
(b) The dg-schemes $RB_nG, n\geq 0$, arrange into a simplicial dg-scheme 
$RB_\bullet G$, containing
$B_\bullet G$ as a closed simplicial subscheme (with trivial dg-structure).
\hfill\break
(c) The simplicial dg-scheme $RBG$ is injective.  

\noindent {\sl Proof:} Part (b) is obvious, part (a) will be proved
in the next subsection.

\vskip .3cm

\noindent {\bf (1.9) $RB_nGL(r)$ and simplicial connections.}
It is easy to understand the meaning of the algebra ${\bf C}[RB_nG]$. 
For any group $G$, points of the
space $\tilde B_nG = \prod_{0\leq i<j\leq n} G$ can be viewed as simplicial 
$G$-connections
on the $n$-simplex $\Delta[n]$, while points of the subscheme $B_nG$ can be  
viewed
as flat connections. 
When $G=GL(r)$ (which assumption we will keep), 
this analogy can be made more precise as follows.

For any associative algebra $R$ let $C^\bullet(\Delta[n], R)$ be the 
(normalized)
simplicial cochain
complex of $\Delta[n]$ with coefficients in $R$. Non-degenerate $p$-faces of 
$\Delta [n]$ are labelled
by sequences $(i_0, ..., i_p)$, $0\leq i_0 < ... < i_p\leq n$, and thus an 
element of
$C^p(\Delta[n], R)$ is a function $\phi$ associating to any such sequence an 
element
$\phi(i_0, ..., i_p)\in R$. The Alexander-Whitney multiplication
$$(\phi\cdot\psi) (i_0, ..., i_{p+q}) = \phi(i_q, ..., i_{p+q}) \psi(i_0, ..., 
i_q),
\quad \phi\in C^p(\Delta[n], R), \psi\in C^q(\Delta[n], R),\leqno (1.9.1)$$
makes $C^\bullet(\Delta[n], R)$ into an associative dg-algebra.  

Now let ${\bf gl}(r)$ be the associative algebra of $r$ by $r$ matrices. A 
point ${\bf g} =
(g_{ij})\in\tilde B_n G$ is nothing but an element of $C^1(\Delta[n], {\bf 
gl}(r))$
whose components are all invertible. The condition for {\bf g} to lie in the
subscheme $B_nG$ can be expressed as
$$d{\bf g} +{\bf  g}\cdot {\bf g} = 0 \quad {\rm in} \quad C^2(\Delta[n], {\bf 
gl}(r)).
\leqno (1.9.2)$$
Further, let $A^\bullet$ be any dg-algebra. An element $\gamma \in
 C^\bullet(\Delta[n], {\bf gl}(r))\otimes A^\bullet$ of degree 1 can be split 
into its components
$\gamma_p\in C^p(\Delta[n], {\bf gl}(r))\otimes A^{1-p}$. Each $\gamma_p$ can 
be
viewed as a collection of matrices $\gamma_{i_0, ..., i_p}$ whose matrix 
elements
belong to $A^{1-p}$. By comparing the formula for the Alexander-Whitney map
 with the definition of the differential in the algebra  ${\bf C}[RB_nG]$ we 
get the
following characterization of  the latter.

\proclaim (1.9.3) Proposition. Let $A^\bullet$ be any dg-algebra. A 
dg-homomorphism
${\bf C}[RB_nG]\to A^\bullet$ is the same as a degree 1 element $\gamma\in
C^\bullet(\Delta[n], {\bf gl}(r))\otimes A^\bullet$ satisfying $d\gamma+ 
\gamma\cdot\gamma = 0$
and such that all the matrices $\gamma_{i_0i_1}$ with entries in $A^0$ are 
invertible.

\vskip .2cm

\noindent {\bf (1.10) Proof of Theorem 1.8.3 (a).} By Proposition 1.2.3, 
it is enough to prove the following fact.

\proclaim (1.10.1) Lemma. Let ${\bf g}\in B_nG \i RB_nG$. Then the tangent 
dg-space
$T_{\bf g} RB_nG$ has no cohomology in degrees other than 0. 

To see the lemma, let us regard {\bf g} as a simplicial local system $V$ on
$\Delta[n]$. Then, we identify  the complex $T_{\bf g} RB_nG$ (concentrated in
non-negative degrees)  with the
shifted and truncated cochain complex
$$C^1(\Delta[n], {\rm End}(V))\to C^2(\Delta[n], {\rm End}(V))\to ...$$
which is clearly exact outside the leftmost term.

\vfill\eject

\centerline {\bf \S 2. The derived space of local systems.}

\vskip 1cm

\noindent {\bf (2.1) Ordinary moduli space.} Let $S$ be a connected
simplicial set and $G$ be a reductive algebraic group, as before.
Consider the scheme $\underline{\rm Hom}(S, BG)$.
It is acted upon
by the group $\prod_{x\in S_0} G$. Two points 
$f,g\in \underline{\rm Hom}(S, BG)$ are equivalent with respect
to this action if and only if the morphisms $f,g$ are elementary
homotopic, i.e., can be obtained as restrictions of one
morphism $F: S\times \Delta[1]\to BG$.
A morphism $f: S\to BG$ is just a rule associating to
any edge $\gamma\in S_1$ an element of $G$ so that for any
2-simplex $\sigma\in S_2$ we have the condition
$g_{\partial_0\sigma}g_{\partial_2\sigma} = g_{\partial_1\sigma}$.
So, geometrically, it can be viewed as a  flat simplicial
$G$-connection on $S$ and the elements of  $\prod_{x\in S_0} G$
are discrete analogs of gauge transformations.

Let $x_0\in S_0$ be a vertex. Define
$$\LS(S, x_0) := \underline{\rm Hom}(S, BG)\biggl/\prod_{x\in S, x\neq x_0} G. 
\leqno (2.1.1)$$
Notice that the action of the subgroup $\prod_{x\neq x_0} G$ is free, so
the quotient exists as an algebraic variety. It is clear that 
$$\LS(S, x_0) = \underline {\rm Hom}(\pi_1(|S|, x_0), G), \leqno (2.1.2)$$
the moduli space of all homomorphisms from the fundamental group
to $G$. This variety can be singular. 

\vskip .3cm

\noindent {\bf (2.2) Derived moduli space.} We keep the notations
and assumptions of the previous subsection.  Choose a $\Delta(G)$-equivariant
injective dg-resolution $RBG$ of $BG$, and consider the dg-scheme
$\underline{\rm Hom}(S, RBG)$. Because of the injectivity, this
is a smooth dg-manifold. The group $\prod_{x\in S_0}G$ acts on this
manifold. 

\proclaim (2.2.1) Proposition. The action of the subgroup $\prod_{x\neq x_0}
G$ on $\underline {\rm Hom}(S, RBG)$ is free.

\noindent {\sl Proof:} The ordinary simplicial scheme underlying $RBG$
 (i.e., the simplicial scheme formed by the spectra of the $0$th graded
components of the dg-algebras ${\bf C}[RB_nG]$ is just $\tilde BG$.
So the ordinary scheme underlying $\underline {\rm Hom}(S, RBG)$
is $\underline{\rm Hom}(S, \tilde BG)$, which is just the product of
as many copies of $G$ as there are nondegenerate 1-simplices in $S$. 
In other words, it the space of all simplicial $G$-connections
on $S$, falt or not. The action of $\prod_{x\in S_0, x\neq x_0}G$
on this space is clearly free. The inclusion of the $0$th
component into any ${\bf Z}_{\leq 0}$-graded dg-algebra is
a morphism of dg-algebras. This means that we have a morphism
(in fact, a fibration) of dg-schemes $\underline{\rm Hom}(S, RBG)
\to\underline {\rm Hom}(S, \tilde BG)$. This morphism is equivariant and
the action on the target is free. So the action on the source is free.

\proclaim (2.2.2) Definition. The derived moduli space of local systems
is defined as
$$R\LS(S, x_0) = \underline{\rm Hom}(S, RBG)\biggl/\prod_{x\in S_0, x\neq x_0} 
G.$$

Clearly, this moduli space is well defined up to a quasiisomorphism. 
Further, a morphism $f: (S, x_0)\to (T, y_0)$ of pointed simplicial sets
induces a morhism $f^*: R\LS_G(T, y_0)\to R\LS(S, x_0)$
of dg-manifolds. 

\proclaim (2.3) Theorem. (a) $R\LS(S, x_0)$ is a smooth dg-manifold, with
$$\pi_0(R\LS(X, x_0)) = \LS(S, x_0).$$
(b) For a {\bf C}-point $[E]$ of $R\LS_G(S, x_0)$ represented by a $G$-local
system $E$ on $S$, the cohomology of the tangent dg-space is found as follows:
$$H^i(T^\bullet_{[E]} R\LS(S, x_0)) = \cases{
H^{i+1}(S, {\rm ad}(E)), i\geq 1\cr
H^1_{res}(S, {\rm ad}(E)), i=0\cr}$$
where $H^1_{res}(S, {\rm ad}(E)) = Z^1/dC^0_{res}$ with $Z^1$ being the
space of 1-cocycles and $C^0_{res}$ being the space of 0-cochains whose
values at $x_0$ is zero. 

\noindent {\sl Proof:} (a) The smoothness follows from the smoothness
of $\underline{\rm Hom}(S, RBG)$ and the freeness of the action. 
The statement about $\pi_0$ follows because $\pi_0$ commutes
with $\underline{\rm Hom}(S, -)$ (exactness of the cokernel with
respect to colimits). 

(b) Let $\nabla$ be a morphism $S\to BG$, i.e., a flat simplicial
connection on $S$, and $E$ be the local system represented by
$\nabla$. We will prove that
$$H^i(T_{[\nabla]}^\bullet \underline{\rm Hom}(S, RBG)) = \cases{
H^{i+1}(S, {\rm ad}(E)), i\geq 1\cr
Z^1(S, {\rm ad}(E)), i=0\cr}. \leqno (2.3.1)$$
This will imply our statement since $R\LS(S, x_0)$ is the quotient
of $\underline{\rm Hom}(S, RBG)$ by the group $\prod_{x\neq x_0}G$.

To prove (2.3.1), we first consider the case when $S=\Delta[n]$
is the $n$-simplex. Then, $\underline{\rm Hom}(\Delta[n], RBG)$
is just $RB_nG$ which is a resolution of $B_nG = 
\underline{\rm Hom}(\Delta[n], BG)$. So the tangent dg-space at
$[\nabla]$ to $\underline{\rm Hom}(\Delta[n], RBG)$ is 
quasiisomorphic to the ordinary tangent space at $[\nabla]$
to $\underline {\rm Hom}(\Delta[n], BG)$ situated in degree 0.

Now, consider the case of general $S$ and represent
$S$ as the union (direct limit) of its simplices and realize accordingly
the dg-scheme $\underline{\rm Hom}(S, RBG)$ as an inverse limit:
$$S = {\rm colim}_{s\in S_n, n\geq 0} \Delta[n], 
\quad \underline{\rm Hom}(S, RBG) = {\rm lim}_{s\in S_n, n\geq 0}
\underline{\rm Hom}(\Delta[n], RBG). \leqno (2.3.2)$$ 
On the level of  tangent dg-spaces this implies:
$$T^\bullet_{[\nabla]} \underline{\rm Hom}(S, RBG) = 
{\rm lim}_{s\in S_n, n\geq 0} T^\bullet_{[\nabla]}\underline{\rm Hom}
(\Delta[n], RBG).\leqno (2.3.3)$$
But because of the injectivity of $RBG$ all the maps in the
diagram whose inverse limit is (2.3.3), are surjective morphisms of
complexes. Therefore the limit is quasiisomorphic to the homotopy inverse 
limit:
$$T^\bullet_{[\nabla]} \underline{\rm Hom}(S, RBG) \sim 
{\rm holim}_{s\in S_n, n\geq 0} T^\bullet_{[\nabla]}\underline{\rm Hom}
(\Delta[n], RBG)\sim {\rm holim}_{s\in S_n, n\geq 0} Z^1(\Delta[n],
{\rm ad}(E)).$$
But the last homotopy inverse limit, if we calculate it via the nerve,
will have exactly the cohomology described in (2.3.1).
Theorem is proved.

\vskip .3cm

\proclaim (2.4) Corollary. A weak equivalence $(S, x_0)\to (T, y_0)$
of pointed simplicial sets induces a quasi-isomorphism of dg-manifolds
$R\LS(T, y_0)\to R\LS(S, x_0)$. 

\noindent {\sl Proof:} This follows from Theorem 2.3 and the
 ``Whitehead theorem"
1.2.3.

\vskip .3cm

\noindent {\bf (2.5) The case of simple local systems.}
Suppose that $E$ is a $G$-local system on $S$ such that ${\rm Aut}(E)=\{1\}$.
If $\nabla\in {\rm Hom}(S, BG)$ is any simplical connection representing $E$,
then the action of the full group $\prod_{x\in S_0}G$ is free
on the neighborhood of $[\nabla]$ in $\underline{\rm Hom}(S, BG)$.
The corresponding formal germ of the quotient is denoted by
${\rm Def}(E)$ and is called the formal deformation space of $E$.
This is a formal scheme with one closed point $[E]$.
In this case the action on a neighborhood of $[\nabla]$
in $\underline{\rm Hom}(S, RBG)$ is also free and we get a smooth
dg-thickening ${\rm RDef}(E)$ which is a formal dg-scheme
with $\pi_0 {\rm RDef}(E) = {\rm Def}(E)$. By factorizing the
equality (2.3.1) by  $\prod_{x\in S_0}G$ we get the following fact. 

\proclaim (2.5.1) Proposition. We have 
$H^i T^\bullet_{[E]} {\rm RDef}(E)) = H^{i+1}(S, {\rm ad}(E))$ for
all $i\geq 0$. Thus the dg-algebra structure on the
ring ${\bf C}[{\rm RDef}(E)]$ makes $T^\bullet_{[E]} {\rm RDef}(E)[-1]
\sim R\Gamma(S, {\rm ad}(E))$
into a weak Lie algebra.

This result provides a ``derived" generalization of the main theorem
of Hinich and Schechtman [H] [HS] [Sch] (for the case of local systems). 
Note that the ring of functions on RDef serves as the cochain
complex of the weak Lie algebra structure on $R\Gamma(S, {\rm ad}(E))$,
so, in particular, the $0$th cohomology of this weak Lie algebra
is the algebra of functions on the ordinary formal
moduli space.

\vskip 2cm

\centerline {\bf References.}

\vskip 1cm

\noindent [F] E.M. Friedlander, Etale Homotopy of Simplicial
Schemes (Ann. Math. Sudies {\bf 104}), Princeton Univ. Press, 1982.

\vskip .2cm

\noindent [H]. V. Hinich, Descent of Deligne groupoids,
preprint alg-geom/9606010

\vskip.2cm

\noindent [HS] V. Hinich, V. Schechtman, Deformation theory and Lie algebra
homology, preprint alg-geom/9405013. 

\vskip .2cm

\noindent [K] M. Kontsevich, Enumeration of rational curves via torus actions, 
in: The Moduli Space of Curves (R. Dijkgraaf, C. Faber, G. van der Geer
Eds.) (Progress in Math. {\bf 129}), p. 335-365, Birkhauser, Boston,
1995. 

\vskip .2cm

\noindent [Man] Y.I. Manin, Gauge Fields and Complex Geometry, Springer-Verlag 
1985.

\vskip .2cm

\noindent [Mun] H.J.Munkholm, DGA algebras as a Quillen model category: 
relation to shm maps, {\it J.Pure Appl. Algebra}, {\bf 13} (1978),
221-232.

\vskip .2cm

\noindent [May] J.P. May, Simplicial Objects in Algebraic Topology, 
Univ. of Chicago Press, 1967. 

\vskip .2cm

\noindent [Q] D. Quillen, On the (co)homology of (co)mmutative rings,
in: Proc. Symp. Pure Math. {\bf XVII} (1970), p. 65-87.

\vskip .2cm

\noindent [Sch] V. Schechtman, Local structure
of moduli spaces, preprint alg-geom/9708008.

\vskip .2cm

\noindent [St] J. Stasheff, Differential graded Lie algebras, quasi-Hopf 
algebras and higher
homotopy algebras, in: Lecture Notes in Math., {\bf 1510}, p. 120-137,
Springer-Verlag, 1992.

\vskip 2cm

{\sl Author's address: Department of Mathematics, Northwestern University, 
Evanston IL 60208,
\hfill\break
email: kapranov@math.nwu.edu}

\bye